\documentclass[twoside,12pt]{article}
\usepackage[]{amsmath,amsthm,amssymb}
\usepackage[latin1]{inputenc}

\begin{document}

\title{{\Large\bf  Discrete  index  transformations with squares of Bessel functions}}

\author{Semyon  YAKUBOVICH}
\maketitle

\markboth{\rm \centerline{ Semyon   YAKUBOVICH}}{}
\markright{\rm \centerline{  DISCRETE TRANSFORMS WITH SQUARES OF BESSEL FUNCTIONS}}

\begin{abstract} {\noindent Discrete analogs of the index  transforms with squares of  Bessel functions of the first and second kind $J_\nu(z),\ Y_\nu(z)$ are introduced and investigated. The corresponding inversion theorems for   suitable classes  of functions and sequences   are established. }

\end{abstract}
\vspace{4mm}



\vspace{4mm}

\section {Introduction and preliminary results}

In this paper we continue to investigate mapping properties of discrete index transforms, introducing the following transformations between suitable sequences $\{a_n\}_{n\ge 1}$ and functions $f$ in terms of the series and integrals which are associated with squares of Bessel functions of the first and second kind  (cf. [1],  Ch. 10), namely,

$$f(x)= \sum_{n=1}^\infty a_n \left[ J^2_{in/2}(x) + Y^2_{in/2}(x) \right],\quad x > 0,\eqno(1.1)$$

$$ a_n= \int_0^\infty \left[ J^2_{in/2}(x) + Y^2_{in/2}(x) \right] f(x) dx,\quad n \in \mathbb{N},\eqno(1.2)$$

$$f(x)= \sum_{n=1}^\infty {a_n\over \cosh(\pi n/2)}  {\rm Re} \left[ J^2_{in/2}(x) \right],\quad x > 0,\eqno(1.3)$$

$$ a_n= {1\over \cosh(\pi n/2)} \int_0^\infty  {\rm Re} \left[ J^2_{in/2}(x) \right] f(x) dx,\quad n \in \mathbb{N},\eqno(1.4)$$

$$f(x)= \sum_{n=1}^\infty {a_n\over \sinh(\pi n/2)} {\rm Im} \left[ J^2_{in/2}(x) \right],\quad x > 0,\eqno(1.5)$$

$$ a_n=  {1\over \sinh(\pi n/2)} \int_0^\infty  {\rm Im} \left[ J^2_{in/2}(x) \right] f(x) dx,\quad n \in \mathbb{N}.\eqno(1.6)$$
 We note   that, for instance,  continuum analogs of the transformations (1.3), (1.4) were considered by the author in [4].  Here  $i$ is the imaginary unit and ${\rm Re},\ {\rm Im}$ denote the real and imaginary parts of a complex-valued function.   Bessel functions $J_\nu(z),\ Y_\nu(z),\ z,\nu \in \mathbb{C}$ of the first and second kind, respectively, are solutions of the Bessel differential equation 

$$  z^2{d^2u\over dz^2}  + z{du\over dz} + (z^2- \nu^2)u = 0.\eqno(1.7)$$
These functions have the following asymptotic behavior at infinity and near the origin
$$ J_\nu(z) = \sqrt{2\over \pi z} \cos \left( z- {\pi\over 4} (2\nu+1)\right)  [1+ O(1/z)], \ z \to \infty,\   |\arg z| <  \pi,\eqno(1.8)$$
$$J_\nu(z) = O( z^{\nu} ), \ z \to 0,\eqno(1.9)$$

$$ Y_\nu(z) = \sqrt{2\over \pi z} \sin \left( z- {\pi\over 4} (2\nu+1)\right)  [1+ O(1/z)], \ z \to \infty,\   |\arg z| <  \pi,\eqno(1.10)$$
$$Y_\nu(z) = O\left( |z|^{-{|{\rm Re}\nu}|} \right), \ z \to 0,\ \nu\neq 0,\eqno(1.11)$$
$$Y_0(z) = O\left(\log(|z|)\right),  \ z \to 0.\eqno(1.12)$$
The sum of squares $J^2_{\nu}(z) + Y^2_{\nu}(z)$ is called the Nicholson function and it is represented by the Nicholson integral (cf. [1], Entry 10.9.30)

$$J^2_{\nu}(z) + Y^2_{\nu}(z) = {8\over \pi^2} \int_0^\infty K_0( 2z\sinh(t)) \cosh(2\nu t) dt, \quad  |\arg z| <  {\pi\over 2},\eqno(1.13)$$
where $K_\nu(z)$ is the modified Bessel function.  Moreover, we mention the following Mellin-Barnes integral representations of the kernels of transformations  (1.1)-(1.6) (see [2], Vol. III, Entries 8.4.19.17, 8.4.19.19, 8.4.20.19, 8.4.20.35)

$$  {2\sqrt\pi\over \cosh(\pi n/2)} {\rm Re} \left[ J^2_{in/2}(x) \right] = {1\over 2\pi i} \int_{\gamma-i\infty}^{\gamma +i\infty} \frac {\Gamma((s+ in)/2)\Gamma((s-in)/2) \Gamma((1-s)/2)}{\Gamma(s/2) \ [\Gamma (1-s/2) ]^2 } x^{-s} ds,\eqno(1.14)$$
where $\ x >0, \ 0 < \gamma < 1, \ n \in \mathbb{N}$ and $\Gamma(z)$ is the Euler gamma function, 

$$  {2\sqrt\pi\over \sinh(\pi n/2)} {\rm Im} \left[ J^2_{in/2}(x) \right] $$

$$=  \sqrt\pi \left[ J_{in/2}(x) Y_{-in/2}(x)  +  J_{-in/2}(x) Y_{in/2}(x) \right] $$
 
$$ = - {1\over 2\pi i} \int_{\gamma-i\infty}^{\gamma +i\infty} \frac {\Gamma((s+ in)/2)\Gamma((s-in)/2) }{\Gamma((s+1)/2) \Gamma (1-s/2) } x^{-s} ds,\eqno(1.15)$$
where $\ x >0, \ 0 < \gamma <  3/2,\ n \in \mathbb{N}$,

$$  { \pi^{5/2} \over \cosh(\pi n/2)} \left[ J^2_{in/2}(x) +  Y^2_{in/2}(x)\right] = {1\over 2\pi i} \int_{\gamma-i\infty}^{\gamma +i\infty} \Gamma\left({s\over 2}\right) \Gamma\left({s+ in\over 2}\right)$$

$$\times  \Gamma\left({s+ in\over 2}\right)\Gamma\left({1-s\over 2}\right) x^{-s} ds,\eqno(1.16)$$
where $\ x >0, \ 0 < \gamma < 1,\ n \in \mathbb{N}$. On the other hand, the Nicholson function can be represented in terms of the modified Bessel function via the integral (see [2], Vol. II, Entry 2.16.3.12)

$$ J^2_{in/2}(x) +  Y^2_{in/2}(x) =   {8 \over \pi^2} \  \cosh\left({\pi n\over 2}\right) \int_0^\infty {K_{in}(t)  dt \over (t^2+ 4x^2)^{1/2}},\quad  x >0.\eqno(1.17)$$
The function $K_{in}(t)$ is the kernel of the discrete Kontorovich-Lebedev transform recently investigated by the author [5].   It can be estimated by virtue of the Lebedev inequality (see [3], p. 219)

$$\left| K_{in}(t) \right| \le A { t^{-1/4}\over \sqrt{ \sinh(\pi n)}},\eqno(1.18)$$
where $ A >0$ is an absolute constant.  Therefore we get from (1.17)

$$ \left| J^2_{in/2}(x) +  Y^2_{in/2}(x)\right| \le  {8 A\over \pi^{2}} \  {\cosh\left(\pi n/2\right) \over  \sqrt{ \sinh(\pi n)}} \int_0^\infty {t^{-1/4}  dt \over (t^2+ 4x^2)^{1/2}}$$

$$=   \left({2 \over \pi^2}\right)^{5/4}   A \  \Gamma\left({3\over 8}\right)\Gamma\left({1\over 8}\right)    x^{-1/4}  \coth^{1/2} \left({\pi n\over 2}\right)$$ 
$$\le  \left({2 \over \pi^2}\right)^{5/4}   A \  \Gamma\left({3\over 8}\right)\Gamma\left({1\over 8}\right)  \coth^{1/2} \left({\pi\over 2}\right)  x^{-1/4} ,$$ 
i.e.  
$$ \left| J^2_{in/2}(x) +  Y^2_{in/2}(x) \right| \le  B  x^{-1/4},\eqno(1.19)$$
where $x >0, B $ is a positive constant.  For the modulus of Bessel function we have the inequality from [4]

$$|J_{i\tau}(x)| \le e^x \sqrt{{\sinh(\pi\tau)\over \pi\tau}},\quad x \ge 0,\  \tau \in \mathbb{R}\backslash \{0\}.\eqno(1.20)$$
Hence we find for $x \ge 0, n \in \mathbb{N}$

$$ { \left| {\rm Re} \left[ J^2_{in/2}(x) \right] \right| \over \cosh(\pi n/2)}  \le  {2 e^{2x}\over \pi n}\   \tanh \left({\pi n\over 2}\right) \le   {2 e^{2x}\over \pi n},\eqno(1.21)$$

$$ { \left| {\rm Im} \left[ J^2_{in/2}(x) \right] \right| \over \sinh (\pi n/2)}  \le  {2 e^{2x}\over \pi n}.\eqno(1.22)$$
Finally, we mention Struve functions ${\bf H}_\nu(z),  {\bf K}_\nu(z)$, which are related with Bessel function of the second kind $Y_\nu(z)$ via the equality (see [1], Entry 11.2.5)

$${\bf H}_\nu(z) =   {\bf K}_\nu(z) + Y_\nu(z).\eqno(1.23)$$
The Struve function  ${\bf K}_\nu(z)$ has the integral representation (cf.  [1], Entry 11.5.2)

$${\bf K}_\nu(z)= {2(z/2)^{\nu}\over \sqrt\pi \Gamma(\nu+1/2)}\int_0^\infty e^{-zt} (1+t^2)^{\nu-1/2} dt,\quad {\rm Re} (z) > 0\eqno(1.24)$$
and behaves at infinity as ${\bf K}_\nu(z) = O(z^{\nu-1}),\ z \to \infty$ (cf. [1], Entry 11.6.1).   The Mellin-Barnes representation for this function is given in [2], Vol. III, Entry 8.4.25.3

$${\bf K}_\nu(x)= {\cos(\pi\nu)\over 4\pi^3 i} \int_{\gamma-i\infty}^{\gamma +i\infty} \Gamma\left({s+\nu\over 2}\right) \Gamma\left({s- \nu\over 2}\right)$$

$$\times  \Gamma\left({s+ \nu+1\over 2}\right)\Gamma\left({1-\nu-s\over 2}\right) \left({x\over 2}\right)^{-s} ds,\quad  |{\rm Re} \nu| < \gamma < 1 -  {\rm Re} \nu. \eqno(1.25)$$
We note that these preliminary results will be applied in the sequel to prove inversion theorems for discrete transformations (1.1)-(1.6).

\section{Inversion theorems} 

We begin with

{\bf Theorem 1}. {\it   Let a sequence $a= \{a_n\}_{n\in \mathbb{N}} \in l_1$,  i.e. 

$$||a||_{l_1}= \sum_{n=1}^\infty  |a_n|  < \infty.\eqno(2.1)$$
Then the discrete transformation $(1.1)$ can be inverted by the formula

$$a_n =  \sinh\left({\pi n\over 2}\right)  \int_0^\infty   \Phi_n(x) f(x) dx,\ n \in \mathbb{N},\eqno(2.2)$$
where the kernel $\Phi_n(x)$ is defined by 

$$\Phi_n(x) = x \int_0^\pi J_0(2 x\cosh(u)) \sinh(2u) \sin(nu) du,\quad x >0,\ n \in \mathbb{N},\eqno(2.3)$$
and  integral  $(2.2)$ converges in the improper sense. }

\begin{proof}  The key ingredient to prove (2.2) will be the following relatively convergent  integral (see [2], Vol. II, Entry 2.13.25.9)

$$ \int_0^\infty x  J_0(2 x\cosh(u)) \left[ J^2_{in/2}(x) + Y^2_{in/2}(x) \right] dx = {2 \sin(nu) \over \pi \sinh(2u) \sinh(\pi n/2)}.\eqno(2.4)$$
In fact, substituting $f$ by formula (1.1) and $\Phi_n(x)$ by (2.3) on the right-hand side of (2.2), we change the order of the proper integration and summation  to obtain

$$\sinh\left({\pi n\over 2}\right) \int_0^\infty   \Phi_n(x) f(x) dx =    \sinh\left({\pi n\over 2}\right) \lim_{T \to \infty} \sum_{m=1}^\infty a_m \int_0^\pi \sinh(2u) \sin(nu)  $$

$$\times \int_0^T x  J_0(2 x\cosh(u)) \left[ J^2_{im/2}(x) + Y^2_{im/2}(x) \right] dx du,\eqno(2.5) $$
where the interchange follows immediately from the assumption (3.1) and the estimate for each fixed $T >0$ by virtue of (1.8), (1.19)

$$\sum_{m=1}^\infty |a_m |\int_0^\pi \sinh(2u) |\sin(nu) | $$

$$\times \int_0^T x  |J_0(2 x\cosh(u))|  \left[ J^2_{im/2}(x) + Y^2_{im/2}(x) \right] dx du$$

$$ \le C  ||a||_{l_1} \int_0^T x^{1/4} dx  \int_0^\pi \sinh(u) \cosh^{1/2}(u) du =  {8\over 15}  C\  T^{5/4} ||a||_{l_1} \left[ \cosh^{3/2}(\pi)-1\right], $$
where $C >0$ is an absolute constant.  The problem is to justify the passage to the limit in (2.5) under the summation sign when $T \to \infty$.  To do this, we appeal to the asymptotic behavior (1.8) of Bessel function $J_0$, the modified Bessel function $K_0$ and integral representation (1.13) of the Nicholson function.  Then for a big enough positive $T$ we have via integration by parts and Entry 2.16.2.2 in [2], Vol. II

$$\left| \int_T^\infty  x  J_0(2 x\cosh(u)) \left[ J^2_{im/2}(x) + Y^2_{im/2}(x) \right] dx\right|$$

$$=  {8\sqrt 2 \over \pi^{5/2}} \left| \int_T^\infty  \sqrt x  \cos \left(  2 x\cosh(u) - {\pi\over 4} \right)  \left[1+ O\left({1\over x}\right)\right] \right.$$

$$\left. \times  \int_0^\infty K_0( 2x\sinh(t)) \cos(m t) dt dx\right|$$

$$\le  {8\sqrt 2 \over \pi^{5/2}} \left|- {\sqrt T\over 2\cosh(u) }  \sin \left(  2 T\cosh(u) - {\pi\over 4} \right)  \right.$$

$$\left. \times  \int_0^\infty K_0( 2T\sinh(t)) \cos(m t) dt\right.$$

$$- {1\over 4\cosh(u)} \int_T^\infty  {1\over \sqrt x}  \sin \left(  2 x\cosh(u) - {\pi\over 4} \right) $$

$$\times  \int_0^\infty K_0( 2x\sinh(t)) \cos(m t) dt dx$$

$$+ {1\over \cosh(u)} \int_T^\infty   \sqrt x \sin \left(  2 x\cosh(u) - {\pi\over 4} \right) $$

$$\left. \times  \int_0^\infty K_1( 2x\sinh(t)) \sinh(t) \cos(m t) dt dx\right|$$

$$+  O \left(\int_T^\infty  {1\over \sqrt x}  \int_0^\infty K_0( 2x t)  dt dx\right) $$

$$\le  {8\sqrt 2 \over \pi^{5/2}} \left[  {1\over 4\cosh(u) \sqrt T}  \int_0^\infty K_0( t)  dt + {1\over 8\cosh(u)} \int_T^\infty  {dx\over x^{3/2}}  \int_0^\infty K_0(  t) dt\right.$$

$$\left.+ {1\over 4 \cosh(u)} \int_T^\infty  {1\over  x^{3/2} } \int_0^\infty K_0( t) \cosh \left({t\over 2x}\right)  dt dx +  O \left({1\over \sqrt T} \right) \right]$$

 $$\le  {8\sqrt 2 \over \pi^{5/2}} \left[  {\pi\over 4\cosh(u) \sqrt T} + {1\over 2 \cosh(u)\sqrt T} \int_0^\infty K_0( t) \cosh \left({t\over 2}\right)  dt+  O \left({1\over \sqrt T} \right) \right]$$

$$=   O \left({1\over \sqrt T} \right),\quad T \to \infty,\eqno(2.6)$$
where the estimate is uniform with respect to  $m \in \mathbb{N}$ and  $u \in [0,\pi].$ Consequently, returning to (2.5), we pass to the limit when $T \to \infty$ under the summation sign and use (2.4) to derive 

$$ \sinh\left({\pi n\over 2}\right) \int_0^\infty   \Phi_n(x) f(x) dx =  {2\over \pi} \ \sinh\left({\pi n\over 2}\right)$$

$$\times  \sum_{m=1}^\infty {a_m \over  \sinh(\pi m/2)} \int_0^\pi  \sin(nu) \sin(mu) du =a_n,$$ 
 completing the proof of Theorem 1. 

\end{proof}

Concerning an inversion formula for the transformation (1.2), we have the following result. 

{\bf Theorem 2}.   {\it Let $f$ be a complex-valued function on $\mathbb{R}_+$ which is represented by the integral 

$$f(x) =  x \int_{-\pi}^\pi J_0(x \cosh(u)) \ \varphi(u)  du,\quad x >0,\eqno(2.7)$$ 
where $ \varphi(u) = \psi(u)\sinh(2u)$ and $\psi$ is a  $2\pi$-periodic function, satisfying the Lipschitz condition on $[-\pi, \pi]$, i.e.

$$\left| \psi(u) - \psi(v)\right| \le C |u-v|, \quad  \forall \  u, v \in  [-\pi, \pi],\eqno(2.8)$$
where $C >0$ is an absolute constant.  Then for all $x >0$ the following inversion formula for  transformation $(1.2)$  holds

$$ f(x)  =    \sum_{n=1}^\infty    \sinh\left({\pi n\over 2}\right)  \Phi_n (x) a_n,\eqno(2.9)$$
where $\Phi_n$ is defined by $(2.3)$.}

\begin{proof}   Plugging the right-hand side of the representation (2.7) in (1.2), we change the order of integration and employ  (2.4) to obtain 

$$a_n  =  {2 \over \pi \sinh(\pi n/2)} \int_{-\pi}^\pi { \varphi(u) \sin(nu)  \over  \sinh(2u)} du.\eqno(2.10)$$
This interchange is valid by virtue of (2.6), which guarantees the uniform convergence of the integral (2.4) with respect to $u \in [-\pi, \pi]$.  Hence, following ideas which are elaborated in [5], 
we substitute $a_n$ by formula (2.10) and $\Phi_n$ by (2.3) into the partial sum of the series (2.9) $S_N(x) $, and  it becomes 
$$S_N(x)  ={x\over \pi}  \sum_{n=1}^N   \int_{-\pi}^\pi J_0(x \cosh(t)) \ \sinh(2t) \sin(nt) dt \int_{-\pi}^\pi   {\varphi(u) \over \sinh(2u)} \sin(nu) du.\eqno(2.11)$$
Hence,  calculating the sum via the known identity 
$$ \sum_{n=1}^N  \sin(nt) \sin(nu) = {1\over 4} \left[  {\sin \left((2N+1) (u-t)/2 \right)\over \sin( (u-t) /2)}  -  {\sin \left((2N+1) (u+t)/2 \right)\over \sin( (u+t) /2)} \right],\eqno(2.12)$$
and invoking the definition of $\varphi$, equality (2.11) turns to be as follows 

$$  S_N(x)  =  {x\over 4 \pi} \   \int_{-\pi}^\pi  J_0(x \cosh(t))  \sinh(2t)   \int_{-\pi}^{\pi}   {\varphi(u)+ \varphi(-u)  \over \sinh(2u)} \  {\sin \left((2N+1) (u-t)/2 \right)\over \sin( (u-t) /2)}  du dt $$

$$=   {x\over 4 \pi} \   \int_{-\pi}^\pi   J_0(x \cosh(t)) \sinh(2t)   \int_{-\pi}^{\pi}  \left[ \psi(u)- \psi(-u) \right]  \  {\sin \left((2N+1) (u-t)/2 \right)\over \sin( (u-t) /2)}  du dt.\eqno(2.13)$$
Since $\psi$ is $2\pi$-periodic, we treat  the latter integral with respect to $u$ as follows 

$$  \int_{-\pi}^{\pi}  \left[ \psi(u)- \psi(-u) \right]  \  {\sin \left((2N+1) (u-t)/2 \right)\over \sin( (u-t) /2)}  du $$

$$=  \int_{ t-\pi}^{t+ \pi}  \left[ \psi(u)- \psi(-u) \right]  \  {\sin \left((2N+1) (u-t)/2 \right)\over \sin( (u-t) /2)}  du $$

$$=  \int_{ -\pi}^{\pi}  \left[ \psi(u+t)- \psi(-u-t) \right]  \  {\sin \left((2N+1) u/2 \right)\over \sin( u /2)}  du. $$
Moreover,

$$ {1\over 2\pi} \int_{ -\pi}^{\pi}  \left[ \psi(u+t)- \psi(-u-t) \right]  \  {\sin \left((2N+1) u/2 \right)\over \sin( u /2)}  du - \left[ \psi(t)- \psi(-t) \right] $$

$$=  {1\over 2\pi} \int_{ -\pi}^{\pi}  \left[ \psi(u+t)- \psi(t) + \psi (-t) - \psi(-u-t) \right]  \  {\sin \left((2N+1) u/2 \right)\over \sin( u /2)}  du.$$
When  $u+t > \pi$ or  $u+t < -\pi$ then we interpret  the value  $\psi(u+t)- \psi(t)$ by  formulas

$$\psi(u+t)- \psi(t) = \psi(u+t-2\pi)- \psi(t - 2\pi),$$ 

$$\psi(u+t)- \psi(t) = \psi(u+t+ 2\pi)- \psi(t +2\pi),$$ 
respectively.  Analogously, the value  $\psi(-u-t)- \psi(-t)$  can be treated.   Then   due to the Lipschitz condition (2.8) we have the uniform estimate
for any $t \in [-\pi,\pi]$

$${\left|  \psi(u+t)- \psi(t) + \psi (-t) - \psi(-u-t) \right| \over | \sin( u /2) |}  \le 2C \left| {u\over \sin( u /2)} \right|.$$
Therefore,  owing to the Riemann-Lebesgue lemma

$$\lim_{N\to \infty } {1\over 2\pi} \int_{ -\pi}^{\pi}  \left[ \psi(u+t)- \psi(-u-t)  - \psi(t) + \psi (-t) \right]  \  {\sin \left((2N+1) u/2 \right)\over \sin( u /2)}  du =  0\eqno(2.14)$$
for all $ t\in [-\pi,\pi].$    Besides, returning to (2.13), we estimate the iterated integral (see (1.8))

$$ \int_{-\pi}^\pi \left|J_0(x \cosh(t))  \sinh(2t) \right|  \int_{ -\pi}^{\pi} \left| \left[ \psi(u+t)- \psi(-u-t)  - \psi(t) + \psi (-t) \right]\right.$$

$$\left.\times   {\sin \left((2N+1) u/2 \right)\over \sin( u /2)}  \right| du dt \le  {4 C_1\over \sqrt x}  \int_{0}^\pi  {\sinh(2t)\over \cosh^{1/2}(t)}    dt   \int_{ -\pi}^{\pi}   \left| {u\over \sin( u /2)} \right| du < \infty,\ x >0,$$
where $C_1 >0$ is constant. Consequently, via  the dominated convergence theorem it is possible to pass to the limit when $N \to \infty$ under the  integral sign, and recalling (2.14), we derive

$$  \lim_{N \to \infty}   {x\over 4 \pi}  \int_{-\pi}^\pi J_0(x \cosh(t))  \sinh(2t)   \int_{ -\pi}^{\pi}  \left[ \psi(u+t)- \psi(-u-t)  - \psi(t) + \psi (-t) \right] $$

$$\times  \  {\sin \left((2N+1) u/2 \right)\over \sin( u /2)}  du dt =  {x\over 4 \pi}  \int_{-\pi}^\pi J_0(x \cosh(t))  \sinh(2t)   $$

$$ \times \lim_{N \to \infty}  \int_{ -\pi}^{\pi}  \left[ \psi(u+t)- \psi(-u-t)  - \psi(t) + \psi (-t) \right]  \  {\sin \left((2N+1) u/2 \right)\over \sin( u /2)}  du dt = 0.$$
Hence, combining with (2.13),  we obtain  by virtue of  the definition of $\varphi$ and $f$

$$ \lim_{N \to \infty}  S_N(x) =   {x\over 2} \   \int_{-\pi}^\pi  J_0(x \cosh(t))  \left[ \varphi (t)+ \varphi(-t) \right] dt = f(x),$$
where the integral (2.7) converges since $\varphi \in C[-\pi,\pi]$.  Thus we established  (2.9), completing the proof of Theorem 2.
 
\end{proof} 

In order to invert discrete transformation (1.3) we will need the following lemmas.

{\bf Lemma 1.}  {\it Let  $u \in [-\pi,\pi],\ n \in \mathbb{N}$. Then the following formula takes place

$$\int_0^\infty x S_{-1,0}\left(2 x\cosh(u)\right)  {\rm Re} \left[ J^2_{in/2}(x) \right] dx =  { \pi \sin(nu) \over  \sinh(2u)  \sinh(\pi n/2) },\eqno(2.15)$$
where $S_{\mu,\nu}(z)$ is the Lommel function (cf. $[1]$,  Sec. $11.9$) and integral $(2.15)$ converges absolutely.}

\begin{proof} The absolute converges of the integral (2.15) follows immediately from (1.8) and asymptotic behavior of the Lommel function at infinity [1], Entry 11.9.9 $S_{\mu,\nu}(z) =
O(z^{\mu-1} ),\\ z \to \infty$. Hence, recalling (1.14), Entry 8.4.27.3 in [2], Vol. III  and Stirling's asymptotic formula for the gamma function [3],  which gives for $0 < \gamma < 1/2$

$$ \frac {\Gamma((s+ in)/2)\Gamma((s-in)/2) \Gamma((1-s)/2)}{\Gamma(s/2) \ [\Gamma (1-s/2)]^2} = O\left( |s|^{\gamma -3/2}\right),\ |s| \to \infty,$$
we appeal to the Parseval equality for the Mellin transform [2], Vol. III to derive 

$$ \int_0^\infty x S_{-1,0}\left(2x\cosh(u)\right)  {\rm Re} \left[ J^2_{in/2}(x) \right] dx$$

$$=  {\cosh(\pi n/2) \over 8\pi^{3/2} i \cosh^2(u)} \int_{\gamma-i\infty}^{\gamma +i\infty}  \Gamma\left(1- {s\over 2}\right)  \Gamma\left({s+ in\over 2}\right)\Gamma\left({s-in\over 2}\right) $$

$$\times \Gamma\left({1-s\over 2}\right)  \left(\cosh(u)\right)^{s} ds.\eqno(2.16)$$
But the latter Mellin-Barnes integral in  (2.16) can be calculated in terms of the Gauss hypergeometric function [3] via  Entry 8.4.49.21 in [2], Vol. III. Indeed, we obtain 

$$ {\cosh(\pi n/2) \over 8\pi^{3/2} i \cosh^2(u)} \int_{\gamma-i\infty}^{\gamma +i\infty}  \Gamma\left(1- {s\over 2}\right)  \Gamma\left({s+ in\over 2}\right)\Gamma\left({s-in\over 2}\right) \Gamma\left({1-s\over 2}\right)  \left(\cosh(u)\right)^{s} ds $$

$$= {\pi n \over 2\cosh(u)\sinh(\pi n/2) }\  {}_2F_1\left({1+in\over 2}, \  {1-in\over 2};\ {3\over 2};\ -\sinh^2(u) \right).\eqno(2.17)$$
Finally, the Gauss hypergeometric function in (2.17) can be simplified due to Entry 7.3.1.91 in [2], Vol. III, and it becomes

$${\pi n \over 2\cosh(u)\sinh(\pi n/2) }  \  {}_2F_1\left({1+in\over 2}, \  {1-in\over 2};\ {3\over 2};\ -\sinh^2(u) \right) =  {\pi  \sin(nu) \over \sinh(2u) \sinh(\pi n/2) }.$$
Thus, combining with (2.16), we end up with (2.15), completing the proof of Lemma 1. 

\end{proof} 

{\bf Lemma 2.}  {\it Let  $x > 0, \ n \in \mathbb{N}$. Then

$${  {\rm Re} \left[ J^2_{in/2}(x) \right] \over \cosh(\pi n/2)}  = {2\over \pi} \int_0^\infty \cos(nt)  {\bf H}_0 \left(2 x\cosh(t)\right) dt,\eqno(2.18)$$
where ${\bf H}_0$ is the Struve  function $(1.23)$ of zero index and integral $(2.18)$ converges absolutely.}

\begin{proof} Since we find from (1.24)

$$ {\bf K}_0 \left( x\right) = {2\over \pi} \int_0^\infty e^{-x\sinh(t)} dt,\eqno(2.19)$$
and $\sqrt x \ Y_0(x),\ x >0$ is bounded (see (1.10), (1.12)),  we have via (1.23) and Entry 2.4.4.4, 2.4.18.4 in [1], Vol. I

$$  \int_0^\infty \left|\cos(nt)  {\bf H}_0 \left( 2x\cosh(t)\right) \right| dt \le    {2\over \pi}   \int_0^\infty   \int_0^\infty  e^{- 2x\cosh(t) \sinh(y)} dy dt $$

$$+ {1\over \sqrt {2x}} \sup_{ u >0} \left| \sqrt u \ Y_0(u) \right|   \int_0^\infty  {dt\over \cosh^{1/2} (t)} =  {2\over \pi}   \int_0^\infty  K_0(2x\sinh(y)) dy  $$

$$+ { \Gamma^2(1/4) \over 4 \sqrt {\pi x} } \sup_{ u >0} \left| \sqrt u \ Y_0(u) \right|  < \infty,\eqno(2.20) $$
where $K_0$ is the modified Bessel function of zero index, we establish the absolute convergence of the integral (2.18).  Then,  appealing to (1.13) and  the following integrals (cf. [2], Vol. II, Entries 2.13.5.2, 2.16.3.12 (cf. (1.17)))

$${4\over \pi} \int_0^\infty \cos(nt) Y_0\left(2x\cosh(t)\right) dt = \left| J_{in}(x) \right|^2-  \left| Y_{in}(x) \right|^2,\eqno(2.21)$$

$$\int_0^\infty K_{in}(2x\sinh(t)) dt = {\pi^2\over 8 \cosh(\pi n/2)} \left[ J^2_{in/2}(x) + Y^2_{in/2}(x) \right],\eqno(2.22) $$
we obtain

$$ {4\over \pi} \int_0^\infty \cos(nt)  {\bf H}_0 \left(2 x\cosh(t)\right) dt = {1 \over  \cosh(\pi n/2)} \left[ J^2_{in/2}(x) + Y^2_{in/2}(x)\right] $$

$$ +  \left| J_{in/2}(x) \right|^2-  \left| Y_{in/2}(x) \right|^2.$$
Hence the final result follows immediately from the relation between Bessel functions of the first and the second kind (see [1], Entry 10.2.3)

$$Y_\nu(z)= {1\over \sin(\pi\nu)}\left[ J_\nu(z)\cos(\pi\nu) - J_{-\nu}(z)\right].\eqno(2.23)$$

\end{proof}

From (2.18), (2.20) we arrive at

{\bf Corollary 1}.  {\it The following inequality holds}

$$ { \left| {\rm Re} \left[ J^2_{in/2}(x) \right] \right|\over \cosh(\pi n/2)} \le  {1\over 2}  \left[  J^2_{0}(x) + Y^2_{0}(x)  +   { \Gamma^2(1/4) \over \pi  \sqrt {\pi x} } \sup_{ u >0} \left| \sqrt u \ Y_0(u) \right|\right].\eqno(2.24)$$

{\bf Theorem 3}. {\it   Let a sequence $a= \{a_n\}_{n\in \mathbb{N}} \in l_1$, i.e. satisfy condition $(2.1)$.  Then the discrete transformation $(1.3)$ can be inverted by the formula

$$a_n =  {1\over \pi^2} \sinh\left(\pi n \right)  \int_0^\infty   \Psi_n(x) f(x) dx,\ n \in \mathbb{N},\eqno(2.25)$$
where the kernel $\Psi_n(x)$ is defined by 

$$\Psi_n(x) = x \int_0^\pi S_{-1,0}(2 x\cosh(u)) \sinh(2u) \sin(nu) du,\quad x >0,\ n \in \mathbb{N},\eqno(2.26)$$
and  integral  $(2.25)$ converges absolutely. }

\begin{proof}  In the same manner as in the proof of Theorem 1 we employ (1.3), (1.8), (1.10), (2.25), (2.26) and Corollary 1 to write the equality

$$ {1\over \pi^2} \sinh\left(\pi n \right)  \int_0^\infty   \Psi_n(x) f(x) dx =  {1\over \pi^2} \sinh\left(\pi n \right) 
 \sum_{m=1}^\infty {a_m\over \cosh(\pi m/2)}  $$
 
 $$\times \int_0^\pi    \sinh(2u) \sin(nu)  \int_0^\infty x   S_{-1,0}(2 x\cosh(u)) {\rm Re} \left[ J^2_{im/2}(x) \right] dx du,\eqno(2.27)$$
where the interchange of the order of integration and summation is permitted by virtue of Fubini's theorem, owing to the estimate
$$\sum_{m=1}^\infty {|a_m|\over \cosh(\pi m/2)}   \int_0^\pi    \sinh(2u) \left|\sin(nu) \right| \int_0^\infty x  \left| S_{-1,0}(2 x\cosh(u)) {\rm Re} \left[ J^2_{im/2}(x) \right] \right|dx du$$

$$\le [\cosh(\pi)-1]  ||a||_{l_1}  \int_0^\infty x  \left| S_{-1,0}(2 x) \right| \left[  J^2_{0}\left({x\over \cosh(\pi)}\right) + Y^2_{0}\left({x\over \cosh(\pi)}\right) \right.$$

$$\left. +   { \Gamma^2(1/4) \cosh(\pi) \over \pi  \sqrt {\pi x} } \sup_{ t >0} \left| \sqrt t \ Y_0(t) \right|\right] dx < \infty.$$
Hence, recalling equality (2.15), we establish (2.25) from (2.27), completing the proof of Theorem 3. 
\end{proof} 

Further, let us consider discrete transformation (1.4). We have

{\bf Theorem 4}.   {\it Let $f$ be a complex-valued function on $\mathbb{R}_+$ which is represented by the integral 

$$f(x) =  x \int_{-\pi}^\pi S_{-1,0}(x \cosh(u)) \ \varphi(u)  du,\quad x >0,\eqno(2.28)$$ 
where $ \varphi(u) = \psi(u)\sinh(2u)$ and $\psi$ is a  $2\pi$-periodic function, satisfying the Lipschitz condition $(2.8)$ on $[-\pi, \pi]$. Then for all $x >0$ the inversion formula for  transformation $(1.4)$  takes place 

$$ f(x)  =    {1\over 2\pi^2} \sum_{n=1}^\infty    \sinh\left(\pi n \right)  \Psi_n (x) a_n,\eqno(2.29)$$
where $\Psi_n$ is defined by $(2.26)$.}

\begin{proof}   Following the same scheme as in the proof of Theorem 2, we substitute  the right-hand side of the representation (2.28) in (1.4), we interchange the order of integration and use  (2.15) to get 

$$a_n  =  { \pi\over  \sinh(\pi n)} \int_{-\pi}^\pi { \varphi(u) \sin(nu)  \over  \sinh(2u)} du.\eqno(2.30)$$
The interchange is allowed  by virtue of Lemma 1, where  the absolute and uniform convergence of the integral (2.15) with respect to $u \in [-\pi, \pi]$ is shown.  Hence we plug
 $a_n$ by formula (2.30) and $\Psi_n$ by (2.26) into the partial sum of the series (2.29) to obtain  
$$S_N(x)  ={x\over \pi}  \sum_{n=1}^N   \int_{-\pi}^\pi S_{-1,0}(x \cosh(t)) \ \sinh(2t) \sin(nt) dt \int_{-\pi}^\pi   {\varphi(u) \over \sinh(2u)} \sin(nu) du.\eqno(2.31)$$
Hence,  making use (2.12) and the definition of $\varphi$, we find 

$$  S_N(x)  =   {x\over 4 \pi} \   \int_{-\pi}^\pi   S_{-1,0}(x \cosh(t)) \sinh(2t)   \int_{-\pi}^{\pi}  \left[ \psi(u)- \psi(-u) \right] $$

$$\times  \  {\sin \left((2N+1) (u-t)/2 \right)\over \sin( (u-t) /2)}  du dt.\eqno(2.32)$$
Hence as in the proof of Theorem 2, we take into account (2.28) and   the definition of $\varphi$ to deduce 

$$ \lim_{N \to \infty}  S_N(x) =   {x\over 2} \   \int_{-\pi}^\pi  S_{-1,0}(x \cosh(t))  \left[ \varphi (t)+ \varphi(-t) \right] dt = f(x).$$
Theorem 4 is proved. 

\end{proof} 

In order to invert discrete transformation (1.5), we will prove the following lemma.

{\bf Lemma 3.}  {\it Let  $u \in [-\pi,\pi],\ n \in \mathbb{N}$. Then the  formula

$$ \int_0^\infty \left[ x {\bf K}_0\left(2 x\cosh(u)\right) - {1\over \pi \cosh(u)}\right] {\rm Im} \left[ J^2_{in/2}(x) \right] dx =  { \sin(nu) \over  \pi \sinh(2u)  \cosh(\pi n/2) }.\eqno(2.33)$$
holds valid, and  integral $(2.33)$ converges absolutely.}

\begin{proof} Indeed, recalling (1.15), we write the representation of the transformation kernel in (1.5), using Entry 2.12.7.4 in [2], Vol. II, in terms of the integral   

$${ {\rm Im} \left[ J^2_{in/2}(x) \right] \over \sinh(\pi n/2)} = - {2\over \pi} \int_0^\infty \cos(nt) J_0( 2x\cosh(t)) dt.\eqno(2.34)$$
Then we have

$${1\over \sinh(\pi n/2)} \int_0^\infty x {\bf K}_0(2x\cosh(u) {\rm Im} \left[ J^2_{in/2}(x) \right] dx =-  {2\over \pi} \lim_{T\to \infty} \int_0^\infty  \cos(nt)$$

$$\times  \int_0^T x {\bf K}_0(2x\cosh(u)) J_0( 2x\cosh(t)) dx dt \eqno(2.35)$$
where the interchange of the order of integration is guaranteed by the estimate (see (1.8), (2.19))

$$  \int_0^\infty \left| \cos(nt)\right| \int_0^T x {\bf K}_0(2x\cosh(u)) \left|J_0( 2x\cosh(t)) \right|dx dt $$

$$\le {1\over \pi \sqrt 2 \cosh(u)}  \sup_{ y >0} \left| \sqrt y   J_0(y) \right|  \int_0^\infty {dt\over \cosh^{1/2}(t)}  \int_0^T {dx\over \sqrt x} $$ 

$$=  {\sqrt T\  \Gamma^2(1/4) \over (2 \pi)^{3/2} \cosh(u)} \  \sup_{ y >0} \left| \sqrt y  J_0(y) \right|. $$
Then using the asymptotic behavior of the Bessel and Struve functions and integration by parts, we motivate the passage to the limit under the integral sign in (2.35) since for a big enough $T$

$$\left| \int_T^\infty x {\bf K}_0(2x\cosh(u) J_0( 2x\cosh(t)) dx\right| = {1\over \pi \sqrt \pi \cosh(u)\cosh^{1/2}(t)} $$

$$\times \left| \int_T^\infty {1\over \sqrt x} \cos\left(2x\cosh(t) -{\pi\over 4}\right) dx + O\left( {1\over \sqrt T}\right)\right| = O\left( {1\over (T\cosh(t))^{1/2}}\right),\ T\to \infty.$$
Meanwhile, the integral with respect to $x$ in (2.35) over $(0,\infty)$ is calculated in [2], Vol. III, Entry 2.7.16.3, and we find

$$\int_0^\infty x {\bf K}_0(2x\cosh(u) J_0( 2x\cosh(t)) dx =  \left[2\pi \cosh(t) (\cosh(t)+\cosh(u))\right]^{-1}.\eqno(2.36)$$
Therefore, combining with (2.35), we obtain (see [2], Vol. I, Entry 2.5.46.6   , Vol. II,  Entry 2.16.6.1)

$$ {1\over \sinh(\pi n/2)} \int_0^\infty x {\bf K}_0(2x\cosh(u) {\rm Im} \left[ J^2_{in/2}(x) \right] dx$$

$$ = -  {1\over \pi^2 }  \int_0^\infty { \cos(nt) dt\over \cosh(t) (\cosh(t)+\cosh(u)) }$$

$$=  {1\over \pi^2 \cosh(u)}  \int_0^\infty  \cos(nt)  \int_0^\infty e^{- y(\cosh(u)+\cosh(t))} dy dt$$

$$ -  {1\over \pi^2 \cosh(u)}  \int_0^\infty  {\cos(nt) \over \cosh(t)} dt $$ 

$$=   {2 \sin(nu) \over \pi  \sinh(\pi n) \sinh(2u)}  -   {1\over  2 \pi \cosh(\pi n/2)\cosh(u)}.\eqno(2.37)$$ 
Finally, from (1.15) and the reciprocal Mellin transform

$${1\over \sinh(\pi n/2)} \int_0^\infty  {\rm Im} \left[ J^2_{in/2}(x) \right] dx = - {1\over 2\cosh(\pi n/2)},$$
and this yields  (2.33). The absolute convergence of the integral follows, in turn,  from (1.8), (1.9) and the asymptotic behavior of the Struve function (cf. [1], Entry 11.6.1) since

$$ x {\bf K}_0\left(2 x\cosh(u)\right) - {1\over \pi \cosh(u)} \sim  - {1\over 4\pi x^2 \cosh^3(u)},\quad x \to \infty.\eqno(2.38)$$
Lemma 3 is proved. 

\end{proof}

Recalling integral representation (2.34), we find  

{\bf Corollary 2}.  {\it The following inequality holds}

$$ { \left| {\rm Im} \left[ J^2_{in/2}(x) \right] \right|\over \sinh(\pi n/2)} \le  {\Gamma^2(1/4) \over 2\pi \sqrt {\pi x}} \sup_{ u >0} \left| \sqrt u \ J_0(u) \right|.\eqno(2.39)$$

 Now we are ready to prove 
 
 {\bf Theorem 5}. {\it   Let a sequence $a= \{a_n\}_{n\in \mathbb{N}} \in l_1$.  Then the discrete transformation $(1.5)$ can be inverted by the formula

$$a_n =  \sinh\left(\pi n \right)  \int_0^\infty   \Omega_n(x) f(x) dx,\ n \in \mathbb{N},\eqno(2.40)$$
where the kernel $\Omega_n(x)$ is defined by 

$$\Omega_n(x) =  \int_0^\pi \left[  x {\bf K}_0\left(2 x\cosh(u)\right) - {1\over \pi \cosh(u)}\right] \sinh(2u) \sin(nu) du,\quad x >0,\ n \in \mathbb{N},\eqno(2.41)$$
and  integral  $(2.40)$ converges absolutely. }

\begin{proof}  In fact, we have from (1.5) and (2.40) 

$$ \sinh\left(\pi n \right)  \int_0^\infty   \Omega_n(x) f(x) dx =   \sinh\left(\pi n \right)  \sum_{m=1}^\infty {a_m\over \sinh(\pi m/2)}  $$
 
 $$\times \int_0^\pi    \sinh(2u) \sin(nu)  \int_0^\infty \left[  x {\bf K}_0\left(2 x\cosh(u)\right) - {1\over \pi \cosh(u)}\right]  {\rm Im} \left[ J^2_{im/2}(x) \right] dx du,\eqno(2.42)$$
where the interchange of the order of integration and summation is valid  by virtue of Fubini's theorem via the asymptotic behavior of the Struve function (2.38), Corollary 2, owing to the estimate
$$\sum_{m=1}^\infty {|a_m|\over \sinh (\pi m/2)}   \int_0^\pi    \sinh(2u) \left|\sin(nu) \right| $$

$$\times \int_0^\infty \left| \left[  x {\bf K}_0\left(2 x\cosh(u)\right) - {1\over \pi \cosh(u)}\right] {\rm Im} \left[ J^2_{im/2}(x) \right] \right|dx du$$

$$\le   {\Gamma^2(1/4) ||a||_{l_1} \over 2\pi^{3/2} } \sup_{ t >0} \left| \sqrt t \ J_0(t) \right| \int_0^\pi    \sinh(2u) $$

$$\times \int_0^\infty \left|   x {\bf K}_0\left(2 x\cosh(u)\right) - {1\over \pi \cosh(u)}\right| {dx du\over \sqrt x} $$

$$=   {2 \Gamma^2(1/4) ||a||_{l_1} \over \pi^{3/2} }  \left[ \cosh^{1/2}(\pi) - 1\right] \sup_{ t >0} \left| \sqrt t \ J_0(t) \right| $$

$$\times \int_0^\infty \left|   x {\bf K}_0\left(2 x\right) - {1\over \pi}\right| {dx \over \sqrt x} < \infty. $$
Consequently, returning to (2.42) e appealing to Lemma 3, we calculate the integral with respect to $x$ by formula (2.33) and then,  due to the orthogonality of trigonometric functions,  end up with the inversion formula (2.40).  Theorem 5 is proved. 

\end{proof} 

Finally,  we will proceed with a proof of  the inversion theorem for transformation (1.6).

{\bf Theorem 6}.   {\it Let $f$ be a complex-valued function on $\mathbb{R}_+$ which is represented by the integral 

$$f(x) =   \int_{-\pi}^\pi \left[ x {\bf K}_0\left(2 x\cosh(u)\right) - {1\over \pi \cosh(u)}\right] \ \varphi(u)  du,\quad x >0,\eqno(2.43)$$ 
where $ \varphi(u) = \psi(u)\sinh(2u)$ and $\psi$ is a  $2\pi$-periodic function, satisfying the Lipschitz condition $(2.8)$ on $[-\pi, \pi]$. Then for all $x >0$ the inversion formula for  transformation $(1.6)$  holds

$$ f(x)  =    {1\over 2} \sum_{n=1}^\infty    \sinh\left(\pi n \right)  \Omega_n (x) a_n,\eqno(2.44)$$
where $\Omega_n$ is defined by $(2.41)$.}

\begin{proof}  Substituting   the right-hand side of  (2.44) in (1.6), we interchange the order of integration and use  (2.33) to derive

$$a_n  =  { 2\over \pi \sinh(\pi n)} \int_{-\pi}^\pi { \varphi(u) \sin(nu)  \over  \sinh(2u)} du.\eqno(2.45)$$
This interchange is verified via the estimate 

$$ {1\over \sinh(\pi n/2)} \int_0^\infty \left| {\rm Im} \left[ J^2_{in/2}(x) \right] \right|\int_{-\pi}^\pi \left| \left[ x {\bf K}_0\left(2 x\cosh(u)\right) - {1\over \pi \cosh(u)}\right] \ \varphi(u) \right| du dx$$

$$\le  {\Gamma^2(1/4) \over 2\pi \sqrt {\pi}} \sup_{ t >0} \left| \sqrt t \ J_0(t) \right| \max_{t \in [-\pi,\pi] } |\varphi(t)|  \int_0^\infty  \left|  x {\bf K}_0\left(2 x\right) - {1\over \pi } \right|  {dx\over \sqrt x} \int_{-\pi}^\pi  \cosh^{-3/2}(u) du$$

$$\le  {\Gamma^2(1/4) \over \sqrt {\pi}} \sup_{ t >0} \left| \sqrt t \ J_0(t) \right| \max_{t \in [-\pi,\pi] } |\varphi(t)|  \int_0^\infty  \left|  x {\bf K}_0\left(2 x\right) - {1\over \pi } \right|  {dx\over \sqrt x} < \infty,$$
which is an immediate consequence of (2.38), (2.39).  Therefore in the same manner as above we find the following representation for the partial sum of the series (2.44) 
$$S_N(x)  ={1\over \pi}  \sum_{n=1}^N   \int_{-\pi}^\pi  \left[ x {\bf K}_0\left(2 x\cosh(t)\right) - {1\over \pi \cosh(t)}\right] \sinh(2t) \sin(nt) dt $$

$$\times \int_{-\pi}^\pi   {\varphi(u) \over \sinh(2u)} \sin(nu) du$$

$$  =   {1\over 4 \pi} \   \int_{-\pi}^\pi    \left[ x {\bf K}_0\left(2 x\cosh(t)\right) - {1\over \pi \cosh(t)}\right] \sinh(2t)   \int_{-\pi}^{\pi}  \left[ \psi(u)- \psi(-u) \right] $$

$$\times  \  {\sin \left((2N+1) (u-t)/2 \right)\over \sin( (u-t) /2)}  du dt.$$
Hence as in the proof of Theorem 2, we take into account (2.43) and   the definition of $\varphi$ to deduce 

$$ \lim_{N \to \infty}  S_N(x) =   {1\over 2} \   \int_{-\pi}^\pi  \left[ x {\bf K}_0\left(2 x\cosh(t)\right) - {1\over \pi \cosh(t)}\right] \left[ \varphi (t)+ \varphi(-t) \right] dt = f(x).$$
Theorem 6 is proved. 

\end{proof}

\bigskip
\centerline{{\bf Acknowledgments}}
\bigskip

\noindent The work was partially supported by CMUP, which is financed by national funds through FCT (Portugal)  under the project with reference UIDB/00144/2020.

\bigskip
\centerline{{\bf References}}
\bigskip
\baselineskip=12pt
\medskip
\begin{enumerate}

\item[{\bf 1.}\ ] NIST Digital Library of Mathematical Functions. http://dlmf.nist.gov/, Release 1.0.17 of 2017-12-22. F. W. J. Olver, A. B. Olde Daalhuis, D. W. Lozier, B. I. Schneider, R. F. Boisvert, C. W. Clark, B. R. Miller and B. V. Saunders, eds.

\item[{\bf 2.}\ ] A.P. Prudnikov, Yu.A. Brychkov and O.I. Marichev, {\it Integrals and Series}. Vol. I: {\it Elementary
Functions}, Vol. II: {\it Special Functions}, Gordon and Breach, New York and London, 1986, Vol. III : {\it More special functions},  Gordon and Breach, New York and London,  1990.

\item[{\bf 3.}\ ] S. Yakubovich, {\it Index Transforms}, World Scientific Publishing Company, Singapore, New Jersey, London and
Hong Kong, 1996.

\item[{\bf 4.}\ ]  S.  Yakubovich, Index transforms with the squares of Bessel functions.  {\it Integral Transforms Spec. Funct. }, {\bf 27}  (2016),  N 12, 981-994.

\item[{\bf 5.}\ ]  S. Yakubovich, Discrete Kontorovich-Lebedev transforms, {\it Ramanujan J.}  DOI 10.1007/s\\ 11139-020-00313-7.

\end{enumerate}

\vspace{5mm}

\noindent S.Yakubovich\\
Department of  Mathematics,\\
Faculty of Sciences,\\
University of Porto,\\
Campo Alegre st., 687\\
4169-007 Porto\\
Portugal\\
E-Mail: syakubov@fc.up.pt\\

\end{document}